\documentclass{ifacconf}

\usepackage{natbib}                    
\usepackage{amsmath,amssymb,amsfonts}  
\usepackage[utf8]{inputenc}				     
\usepackage{lmodern}
\usepackage{enumitem}					         
\usepackage{colonequals}               
\usepackage{textcomp}
\usepackage[T1]{fontenc}               
\usepackage{microtype}
\usepackage{dsfont}
\usepackage{graphicx}

\newcounter{thmcnt}

\newcounter{assumptioncnt}
\newtheorem{theorem}[thmcnt]{Theorem}

\newtheorem{assumption}[assumptioncnt]{Assumption}

\makeatletter
\def\theorem@headerfont{\upshape\bfseries}
\makeatother

\renewcommand{\u}{\mathbf{u}}
\renewcommand{\v}{\mathbf{v}}

\newcommand{\x}{\mathbf{x}}
\newcommand{\y}{\mathbf{y}}

\newcommand{\W}{\mathcal{W}}
\renewcommand{\L}{\mathcal{L}}
\renewcommand{\O}{\mathcal{O}}
\newcommand{\real}{\mathbb{R}}
\newcommand{\1}{\mathds{1}}
\newcommand{\df}{\nabla\! f}
\newcommand{\Df}{\bar\nabla\! f}
\newcommand{\tp}{\mathsf{T}}
\newcommand{\trace}{\mathrm{trace}}
\newcommand{\defeq}{\colonequals}
\newcommand{\bmat}[1]{\begin{bmatrix}#1\end{bmatrix}}
\newcommand{\avg}{\mathsf{avg}}
\newcommand{\dis}{\mathsf{dis}}
\newcommand{\cond}{\mathsf{cond}}
\renewcommand{\forall}{\text{for all }}

\DeclareMathOperator*{\minimize}{minimize}

\usepackage{algorithm,algpseudocode,xpatch}
\renewcommand{\thealgorithm}{}
\def\mathllap{\mathpalette\mathllapinternal}
\def\mathllapinternal#1#2{\llap{$\mathsurround=0pt#1{#2}$}}
\algrenewcommand\algorithmicindent{0.8em}  
\makeatletter
\xpatchcmd{\algorithmic}{\itemsep\z@}{\itemsep=0.5ex plus2pt}{}{}
\makeatother

\begin{document}

\begin{frontmatter}

\title{A Distributed Optimization Algorithm over Time-Varying Graphs with Efficient Gradient Evaluations\thanksref{footnoteinfo}}

\thanks[footnoteinfo]{This material is based upon work supported by the National Science Foundation under Grant No. 1656951 and 1750162.}

\author[1]{Bryan Van Scoy}
\author[1,2]{Laurent Lessard}

\address[1]{Wisconsin Institute for Discovery}
\address[2]{Department of Electrical Engineering\\
  University of Wisconsin--Madison, Madison, WI~53706, USA
  \texttt{\{vanscoy,laurent.lessard\}@wisc.edu}}

\begin{abstract}
We propose an algorithm for distributed optimization over time-varying communication networks. Our algorithm uses an optimized ratio between the number of rounds of communication and gradient evaluations to achieve fast convergence. The iterates converge to the global optimizer at the same rate as centralized gradient descent when measured in terms of the number of gradient evaluations while using the minimum number of communications to do so. Furthermore, the iterates converge at a near-optimal rate when measured in terms of the number of communication rounds. We compare our algorithm with several other known algorithms on a distributed target localization problem.
\end{abstract}

\end{frontmatter}

\section{Introduction}

We consider the distributed optimization problem
\begin{align}\label{eq:distropt}
\minimize_{x\in\real^d} \ f(x) \quad\text{where}\quad f(x) \defeq \frac{1}{n} \sum_{i=1}^n f_i(x).
\end{align}
Associated with each agent $i\in\{1,\ldots,n\}$ is the local objective function $f_i : \real^d\to\real$ where $n$ is the number of agents and $d$ is the dimension of the problem. The goal is for all the agents to calculate the global optimizer using only local communications and gradient evaluations.

Many algorithms have been proposed recently to solve the distributed optimization problem. Some examples include distributed gradient descent by~\cite{DGD}, EXTRA by \cite{EXTRA}, AugDGM by~\cite{AugDGM}, NIDS by~\cite{NIDS}, DIGing by~\cite{DIGing} and~\cite{QuLi}, Exact Diffusion by~\cite{ExactDiffusion1}, and SVL by~\cite{distralg} among others. In each algorithm, agents do the following at each step:
\begin{itemize}
\item communicate state variables with local neighbors,
\item evaluate the local gradient $\df_i$, and
\item update local state variables.
\end{itemize}
Each algorithm alternates between these three steps and therefore uses the same number of communications and local gradient evaluations. In this paper, however, we allow this ratio to depend on the properties of the objective function and the communication network. To characterize the convergence properties of our algorithm, we use the following notions of time.
\begin{itemize}[itemsep=1mm]
\item We define a \emph{step} as one round of communication and at most one gradient evaluation.
\item We define an \emph{iteration} as $m$ rounds of communication and one gradient evaluation.
\end{itemize}
In other words, an iteration consists of $m$ steps where each step is at least as simple as that of the algorithms previously mentioned. We assume that local state updates have neglible cost and can therefore be performed any number of times per step or iteration.

For example, consider an algorithm that updates as follows:\\[2mm]
{\setlength{\tabcolsep}{1.6mm}
\begin{tabular}{r|ccc|ccc|ccc|}
iteration & 1 & & & 2 & & & 3 & & \\
step & 1 & 2 & 3 & 4 & 5 & 6 & 7 & 8 & 9 \\ \hline
communication & \checkmark & \checkmark & \checkmark & \checkmark & \checkmark & \checkmark & \checkmark & \checkmark & \checkmark \rule{0pt}{2.2ex} \\
gradient evaluation & & & \checkmark & & & \checkmark & & & \checkmark
\end{tabular}}

This algorithm performs three rounds of communication per gradient evaluation, so $m=3$.

\textbf{Main contributions.} In this work, we propose a novel decentralized algorithm for solving~\eqref{eq:distropt}. Instead of using the same number of communication rounds as gradient evaluations, our algorithm sets the ratio between these using global problem parameters. We show the following:
\begin{itemize}[itemsep=2mm]
\item[(1)] The iterates of our algorithm converge to the optimizer with the same rate as centralized gradient descent in terms of number of the iterations. Furthermore, our algorithm achieves this using the minimum number $m$ of communications per gradient evaluation.
\item[(2)] The iterates of our algorithm converge to the optimizer with a near-optimal rate in terms of the number of steps, despite not evaluating the gradient at each step.
\end{itemize}
A decentralized algorithm can trivially obtain the same rate as centralized gradient descent if we use an infinite number of communication rounds per iteration (i.e., $m\to\infty$) since then every agent can compute an exact average at each iteration (and therefore can evaluate the global gradient). We show, however, that our algorithm achieves the same rate with a \emph{finite} number of communication rounds per iteration, and we characterize precisely how many communication rounds are required.

To prove convergence of our algorithm, we make the following assumptions.
\begin{itemize}[itemsep=2mm]
\item The local objective functions satisfy a contraction property that is weaker than assuming smoothness and strong convexity.
\item The communication network may be time-varying and either directed or undirected as long as it is sufficiently connected and the associated weight matrix is doubly stochastic at each step.
\end{itemize}

Perhaps the algorithm most similar to ours is the multi-step dual accelerated (MSDA) algorithm by \cite{scaman2017}. This algorithm also adjusts the ratio between the number of communication rounds and gradient evaluations to achieve fast convergence. The MSDA algorithm is provably optimal in terms of both the number of communications and gradient evaluations when the objective function is smooth and strongly convex and the communication network is fixed. Compared to our algorithm, the MSDA algorithm achieves an accelerated rate of convergence by making stronger assumptions on both the objective function and the communication network while we prove a non-accelerated rate using weaker assumptions.

The remainder of the paper is organized as follows. We first set up the distributed optimization problem along with our assumptions in Section~\ref{sec:setup}, and then present our algorithm along with its main convergence result in Section~\ref{sec:main_results}. We then compare our algorithm with several others on a distributed target localization problem in Section~\ref{sec:simulations}, and conclude in Section~\ref{sec:conclusion}. To simplify the presentation, we defer the main convergence proof to Appendix~\ref{sec:proof_theorem1}.

\textbf{Notation.} We use subscript $i$ to denote the agent and superscript $k$ to denote the iteration. We denote the all-ones vector by $\1\in\real^n$ and the identity matrix by $I_n\in\real^{n\times n}$. We use $\|\cdot\|$ to denote the $2$-norm of a vector as well as the induced $2$-norm of a matrix.

\section{Problem setup}\label{sec:setup}

We now discuss the assumptions on the objective function and the communication network that we make in order to solve the distributed optimization problem~\eqref{eq:distropt}.

\subsection{Objective function}

\begin{assumption}\label{assumption:functions}
The distributed optimization problem~\eqref{eq:distropt} has an optimizer $x^\star\in\real^d$. Furthermore, there exists a \emph{stepsize} $\alpha>0$ and \emph{contraction factor} $\rho\in(0,1)$ such that
\begin{align}\label{eq:T_contraction}
\|x-x^\star-\alpha\,(\df_i(x)-\df_i(x^\star))\|\le\rho\,\|x-x^\star\|
\end{align}
for all $x\in\real^d$ and all $i\in\{1,\ldots,n\}$.
\end{assumption}

Each $\df_i(x^\star)$ is in general nonzero, although we have
\begin{align}\label{eq:gradzero}
\sum_{i=1}^n \df_i(x^\star) = 0.
\end{align}
Assumption~\ref{assumption:functions} also implies that
\begin{align*}
\|x-x^\star-\alpha\,\df(x)\|\le\rho\,\|x-x^\star\| \quad\forall x\in\real^d,
\end{align*}
so the global objective function satisfies the same property as the local functions. Assumption~\ref{assumption:functions} holds if the local functions satisfy a one-point smooth and strong convexity property as described in the following~proposition.

\begin{prop}\label{prop:one_point_convex}
Let $0<\mu\le L$, and suppose each local function $f_i$ is one-point $\mu$-smooth and $L$-strongly convex with respect to the global optimizer, in other words,
\begin{align*}
\mu\|x-x^\star\|^2 \le \bigl(\df_i(x)-\df_i(x^\star)\bigr)^\tp (x-x^\star) \le L\|x-x^\star\|^2
\end{align*}
for all $x\in\real^d$ and all $i\in\{1,\ldots,n\}$. Then~\eqref{eq:T_contraction} holds with stepsize $\alpha=\tfrac{2}{L+\mu}$ and contraction factor $\rho=\tfrac{L-\mu}{L+\mu}$.
\end{prop}

Assumption~\ref{assumption:functions} also holds under the stronger assumption that each $f_i$ is $\mu$-smooth and $L$-strongly convex, meaning that
\begin{align*}
\mu\|x-y\|^2 \le \bigl(\df_i(x)-\df_i(y)\bigr)^\tp (x-y) \le L\|x-y\|^2
\end{align*}
for all $x,y\in\real^d$ and all $i\in\{1,\ldots,n\}$.

\subsection{Communication network}

To characterize the communication among agents, we use a gossip matrix defined as follows.\medskip

\begin{defn}[Gossip matrix]
We say that the matrix $W = \{w_{ij}\} \in \real^{n\times n}$ is a \emph{gossip matrix} if $w_{ij}=0$ whenever agent $i$ does not receive information from agent~$j$. We define the \emph{spectral gap} $\sigma\in\real$ of a gossip matrix $W$ as
\begin{align}\label{eq:gap}
\sigma \defeq \|W - \tfrac{1}{n} \1\1^\tp\|.
\end{align}
Furthermore, we say that $W$ is \emph{doubly-stochastic} if both $W\1=\1$ and $\1^\tp W=\1^\tp$.
\end{defn}

The spectral gap characterizes the connectivity of the communication network. In particular, a small spectral gap corresponds to a well-connected network and vice versa. One way to obtain a gossip matrix is to set $W = I-\L$ where $\L$ is the (possibly weighted) graph Laplacian. We make the following assumption about the gossip matrix.\medskip

\begin{assumption}[Communication network]\label{assumption:W}
There exists a scalar $\sigma\in(0,1)$ such that each agent $i\in\{1,\ldots,n\}$ has access to the $i\textsuperscript{th}$ row of a doubly-stochastic gossip matrix $W$ with spectral gap at most $\sigma$ at each step of the algorithm.
\end{assumption}

Time-varying communication networks that are either directed or undirected can satisfy Assumption~\ref{assumption:W} as long as the associated gossip matrix is doubly stochastic with a known upper bound on its spectral gap. See~\cite{FDLA} for how to optimize the weights of the gossip matrix to minimize the spectral gap, and see~\cite{nedic2015} for distributed optimization over non-doubly-stochastic networks using the push-sum protocol.

\subsection{Centralized gradient descent}

The (centralized) gradient descent iterations are given by
\begin{align}\label{eq:fixed_point_iteration}
x^{k+1} = x^k - \alpha\,\df(x^k)
\end{align}
where $\alpha>0$ is the stepsize. Under Assumption~\ref{assumption:functions}, this method converges to the optimizer linearly with rate $\rho$. In other words, $\|x^k-x^\star\| = \O(\rho^k)$. While this method could be approximated in a decentralized manner using a large number of steps per iteration (so that every agent can compute the average gradient at each iteration), we show that our algorithm achieves the same convergence rate using the minimal number $m$ of necessary rounds of communication per gradient evaluation.

\section{Main Results}\label{sec:main_results}

To solve the distributed optimization problem, we now introduce our algorithm, which depends on the stepsize $\alpha$, contraction factor $\rho$, and spectral gap $\sigma$.

\begin{algorithm}[H]
\floatname{algorithm}{Algorithm}
\small\caption{}
\begin{algorithmic}
\State\textbf{Parameters:} stepsize $\alpha>0$, contraction factor $\rho\in(0,1)$, and spectral gap $\sigma\in(0,1)$.
\State\textbf{Inputs:} local functions $f_i : \real^d\to\real$ on agent $i\in\{1,\ldots,n\}$, gossip matrices $\{w_{ij}^{k\ell}\}$ at iteration $k$ and communication round $\ell$.
\State\textbf{Initialization:}
\begin{itemize}[itemsep=3pt]
\item Each agent $i\in\{1,\ldots,n\}$ chooses $x_i^0,y_i^0 \in \real^d$ such that $\sum_{i=1}^n y_i^0 = 0$ (for example, $y_i^0 = 0$).
\item Define the number of communications per iteration
\[
m \defeq \minimize_{r\ge\rho,\,s\ge\sigma} \ \bigl\lceil \log_s\bigl(\tfrac{\sqrt{1+r}-\sqrt{1-r}}{2}\bigr) \bigr\rceil.
\]
\smallskip
\end{itemize}
\For {iteration $k=0,1,2,\ldots$}
  \For {agent $i\in\{1,\ldots,n\}$}
    \State $v_{i,0}^k = x_i^k$
		\For {step $\ell=1,\ldots,m$}
		  \State $v_{i,\ell}^k = \sum_{j=1}^n w_{ij}^{k\ell}\,v_{j,\ell-1}^k \hfill \text{(local communication)}$
    \EndFor
    \State $\hphantom{x_i^{k+1}}\mathllap{u_i^k} = v_{i,m}^k - \alpha\,\df_i(v_{i,m}^k) \hfill \text{(local gradient evaluation)}$
	  \State $\hphantom{x_i^{k+1}}\mathllap{y_i^{k+1}} = y_i^k + x_i^k - v_{i,m}^k \hfill \text{(local state update)}$
	  \State $\hphantom{x_i^{k+1}}\mathllap{x_i^{k+1}} = u_i^k - \sqrt{1-\rho^2}\,y_i^{k+1} \hfill \text{(local state update)}$
  \EndFor
\EndFor\\
\Return $x_i^k\in\real^d$ is the estimate of $x^\star$ on agent $i$ at iteration $k$
\end{algorithmic}
\end{algorithm}

At iteration $k$ of the algorithm, agent $i$ first communicates with its local neighbors $m$ times using the gossip matrices $\{W^{k,\ell}\}_{\ell=1}^m$, then evaluates its local gradient $\df_i$ at the point resulting from the communication, and finally updates its local state variables $x_i^k$ and $y_i^k$. The output of the algorithm is $x_i^k$, which is the estimate of the optimizer $x^\star$ of the global objective function~$f$. Note that agents are required to know the global parameters $\rho$ and~$\sigma$ so that they can calculate the number of communication rounds $m$.

\begin{figure}
\includegraphics[width=\columnwidth]{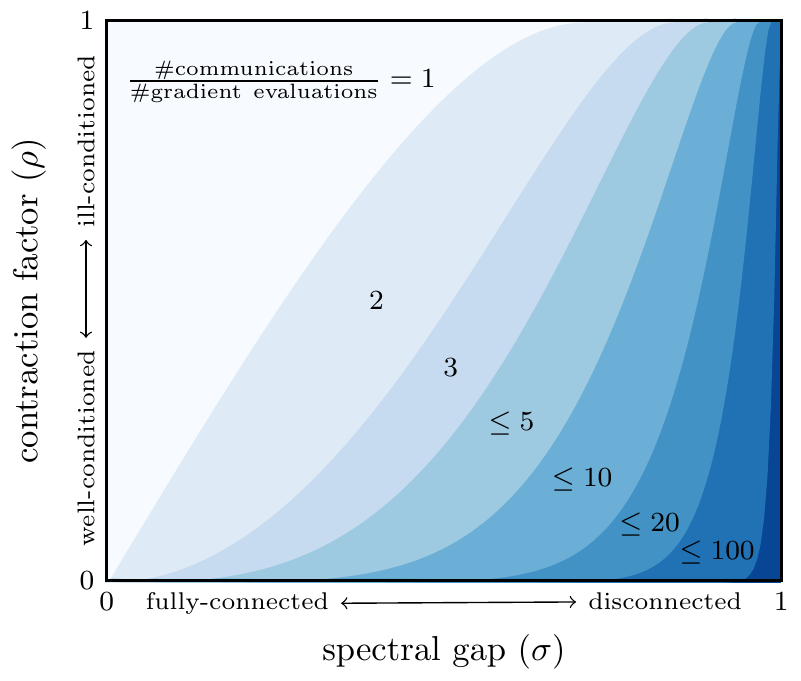}
\caption{Ratio between the number of communications and gradient evaluations as a function of the spectral gap~$\sigma$ and the contraction factor $\rho$. The color indicates the ratio from light (small ratio) to dark (large ratio).\bigskip}
\label{fig:communication}
\end{figure}

For a given contraction factor $\rho$ and spectral gap $\sigma$, agents perform $m$ consecutive rounds of communication at each iteration where
\begin{align}\label{eq:m}
m \defeq \minimize_{r\ge\rho,\,s\ge\sigma} \ \bigl\lceil \log_s\bigl(\tfrac{\sqrt{1+r}-\sqrt{1-r}}{2}\bigr) \bigr\rceil.
\end{align}
This is the minimum integer number of communication rounds so that the spectral gap of the $m$-step gossip matrix $\prod_{\ell=1}^m W^{k,\ell}$ at iteration $k$ is no greater than $\tfrac{\sqrt{1+\rho} - \sqrt{1-\rho}}{2}$. Since only one gradient evaluation is performed per iteration, this adjusts the ratio between the number of communications and gradient evaluations as shown in Figure~\ref{fig:communication}. In particular, the algorithm uses a single communication per gradient evaluation when the network is sufficiently connected ($\sigma$ small) and the objective function is ill-conditioned ($\rho$ large). As the network becomes more disconnected and/or the objective function becomes more well-conditioned, the algorithm uses more communications per gradient evaluation in order to keep the ratio at the optimal operating point.

We now present our main result, which states that the iterates of each agent converge to the global optimizer linearly with a rate equal to the contraction factor $\rho$. We prove the result in Appendix~\ref{sec:proof_theorem1}.\medskip

\begin{figure*}
\includegraphics[width=0.98\columnwidth]{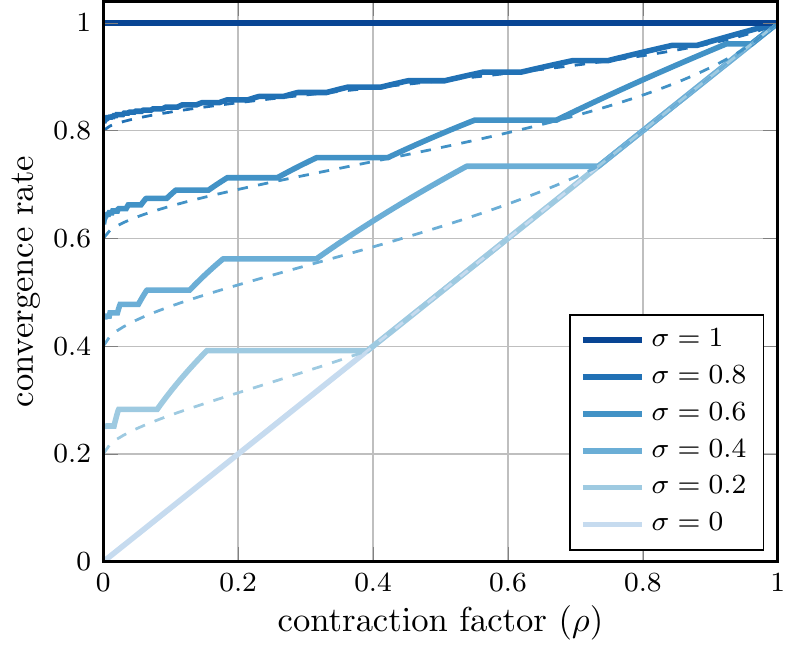}\hfill%
\includegraphics[width=0.98\columnwidth]{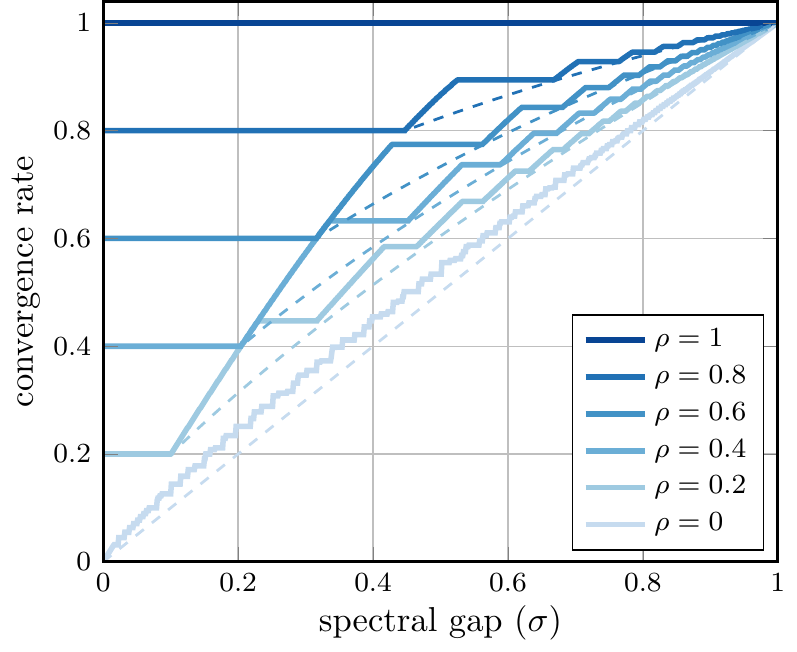}
\caption{Convergence rate in terms of the number of steps as a function of the contraction factor $\rho$ and spectral gap $\sigma$. Solid lines indicate our algorithm while dashed lines indicate the optimal algorithm SVL by~\cite{distralg}. Since our algorithm converges with rate $\rho$ with respect to the number of iterations and performs $m$ steps per iteration, the convergence rate with respect to the number of steps is $\rho^{1/m}$ where $m$ is defined in~\eqref{eq:m}.\bigskip}
\label{fig:contractionFactor}
\end{figure*}

\begin{theorem}[Main result]\label{thm:convergence_rate}
Suppose Assumptions~\ref{assumption:functions} and~\ref{assumption:W} hold for some point $x^\star\in\real^d$, stepsize $\alpha>0$, contraction factor $\rho\in(0,1)$, and spectral gap $\sigma\in(0,1)$. Then the iterate sequence $\{x_i^k\}_{k\ge 0}$ of each agent $i\in\{1,\ldots,n\}$ in our algorithm converges to the optimizer $x^\star$ linearly with rate $\rho$. In other words,
\begin{align}
\|x_i^k-x^\star\| = \O(\rho^k) \quad\forall i\in\{1,\ldots,n\}.
\end{align}
\end{theorem}

Theorem~\ref{thm:convergence_rate} states that the iterates of our algorithm converge to the optimal solution of~\eqref{eq:distropt} in a decentralized manner at the \emph{same} rate as centralized gradient descent~\eqref{eq:fixed_point_iteration} in terms of the number of iterations. In other words, the algorithm converges just as fast (in the worst case) as if each agent had access to the information of all other agents at every iteration. Instead of communicating all this information, however, it is sufficient to only perform $m$ rounds of communication where $m$ is defined in~\eqref{eq:m}.

The convergence rate in Theorem~\ref{thm:convergence_rate} is in terms of the number of iterations. To compare the performance of our algorithm in terms of the number of steps, we plot the convergence rate per step in Figure~\ref{fig:contractionFactor}. For comparison, we also plot the rate of the algorithm SVL by~\cite{distralg}. This algorithm is designed to optimize the convergence rate per step and requires agents to compute their local gradient at each step of the algorithm. In contrast, our algorithm is slightly slower than the optimal algorithm but uses far fewer computations since local gradients are only evaluated once every $m$ steps.

\section{Application: Target Localization}\label{sec:simulations}

\begin{figure}
\includegraphics[width=\columnwidth]{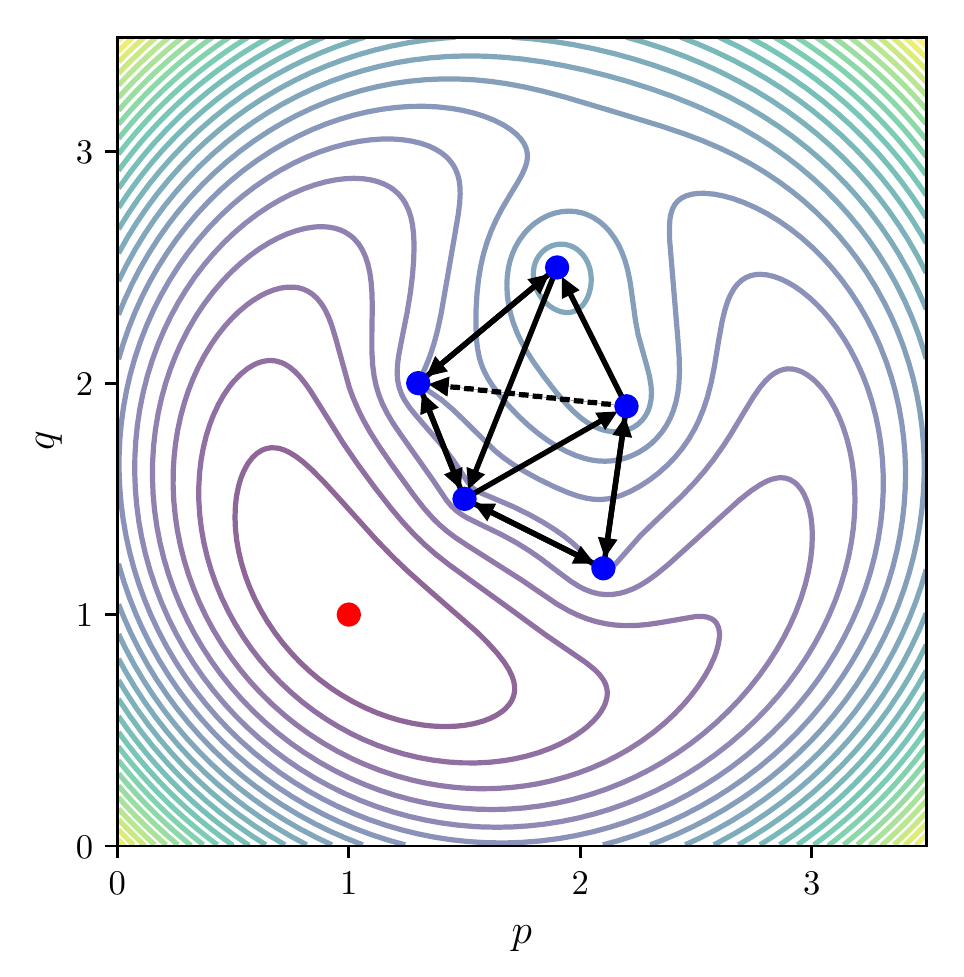}
\caption{Setup of the target localization problem. The position $(p_i,q_i)\in\real^2$ of agent $i\in\{1,\ldots,5\}$ is denoted by a blue dot with the position of the target in red at $(p^\star,q^\star)=(1,1)$. The black arrows indicate the flow of information with an arrow from agent $i$ to $j$ if agent~$j$ receives information from agent $i$. The dashed arrow indicates the link that varies in time. The smooth curves are the contour lines of the objective function for the distributed nonlinear least-squares problem in~\eqref{eq:DNLS}. Note that the problem is nonconvex since the level sets are nonconvex.}
\label{fig:target_localization}
\end{figure}

To illustrate our results, we use our algorithm to solve the distributed target localization problem illustrated in Figure~\ref{fig:target_localization}, which is inspired by the example in Section 18.3 of the book by~\cite{linalg}. We assume each agent (blue dot) can measure its distance (but not angle) to the target (red dot) and can communicate with local neighbors.

Suppose agents are located in a two-dimensional plane where the location of agent $i\in\{1,\ldots,n\}$ is given by $(p_i,q_i)\in\real^2$. Each agent knows its own position but \emph{not} the location of the target, denoted by $x^\star=(p^\star,q^\star)\in\real^2$. Agent $i$ is capable of measuring its distance to the target,
\begin{align*}
r_i = \sqrt{(p_i-p^\star)^2 + (q_i-q^\star)^2}.
\end{align*}
The objective function $f_i : \real^2\to\real$ associated to agent $i$ is
\begin{align*}
f_i(p,q) = \tfrac{1}{2}\,\bigl(\sqrt{(p_i-p)^2 + (q_i-q)^2} - r_i\bigr)^2.
\end{align*}
Then in order to locate the target, the agents cooperate to solve the distributed nonlinear least-squares problem
\begin{align}\label{eq:DNLS}
\minimize_{p,q\in\real} \ \frac{1}{n} \sum_{i=1}^n f_i(p,q).
\end{align}
Agents can communicate with local neighbors as shown in Figure~\ref{fig:target_localization}. To simulate randomly dropped packets from agent $4$ to agent $1$, the gossip matrix at each iteration is randomly chosen from the set
\begin{align*}
W \in \left\{
  \bmat{
             0~ & \tfrac{3}{8}~ & \tfrac{1}{4}~ &            0~ & \tfrac{3}{8} \\[1mm]
  \tfrac{1}{8}~ &            0~ & \tfrac{3}{4}~ & \tfrac{1}{8}~ &            0 \\[1mm]
             0~ & \tfrac{5}{8}~ &            0~ & \tfrac{3}{8}~ &            0 \\[1mm]
  \tfrac{3}{8}~ &            0~ &            0~ &            0~ & \tfrac{5}{8} \\[1mm]
  \tfrac{1}{2}~ &            0~ &            0~ & \tfrac{1}{2}~ &            0}, \
  \bmat{
             0~ & \tfrac{1}{2}~ & \tfrac{1}{4}~ &            0~ & \tfrac{1}{4} \\[1mm]
  \tfrac{1}{4}~ &            0~ & \tfrac{3}{4}~ &            0~ &            0 \\[1mm]
             0~ & \tfrac{1}{2}~ &            0~ & \tfrac{1}{2}~ &            0 \\[1mm]
  \tfrac{1}{4}~ &            0~ &            0~ &            0~ & \tfrac{3}{4} \\[1mm]
  \tfrac{1}{2}~ &            0~ &            0~ & \tfrac{1}{2}~ &            0}\right\}.
\end{align*}
Both gossip matrices satisfy Assumption~\ref{assumption:W} with maximum spectral gap $\sigma \approx 0.7853$.

\begin{figure}
\includegraphics[width=\columnwidth]{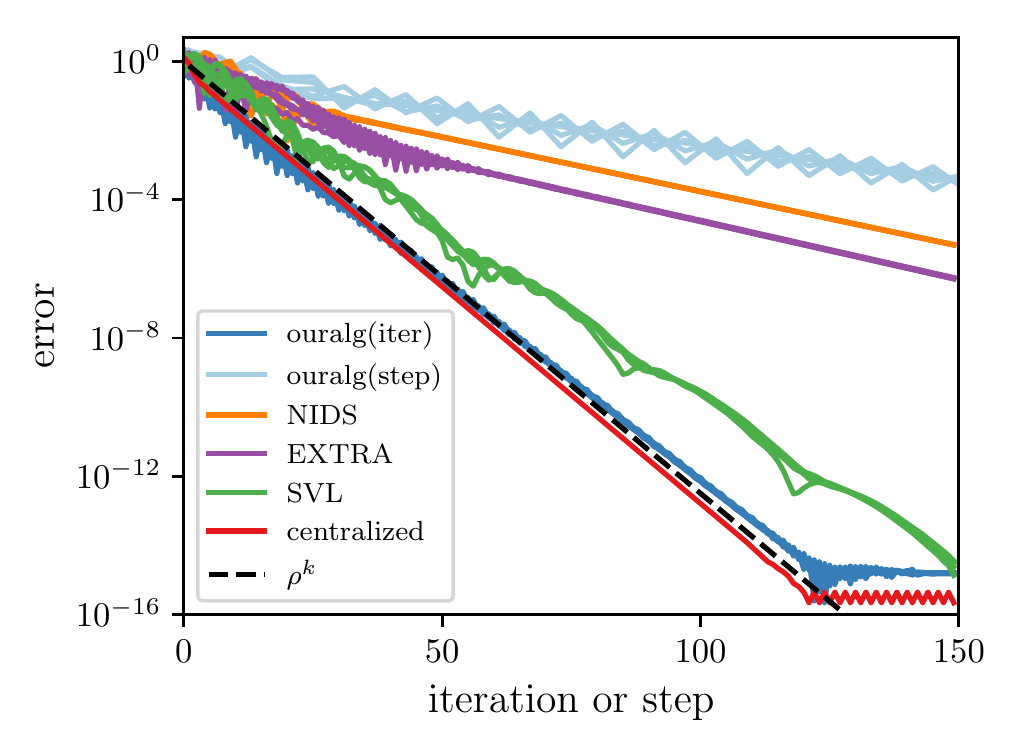}
\caption{Plot of the error for the target localization problem. The blue lines indicate the error $\|x_i^k-x^\star\|$ for each of the five agents computed using our proposed decentralized algorithm as a function of either the iteration (dark blue) or step (light blue); iterations and steps are equivalent for each of the other algorithms. The red line indicates the error using centralized gradient descent~\eqref{eq:fixed_point_iteration}. Our algorithm performs one gradient evaluation and $m=6$ communications per iteration.\smallskip}
\label{fig:target_localization_error}
\end{figure}

We choose the stepsize to optimize the asymptotic rate of convergence. In particular, the estimate of each agent becomes arbitrarily close to the target as $k\to\infty$, so the optimal stepsize is $\alpha=\tfrac{2}{\lambda_1+\lambda_2}$ where $\lambda_1$ and $\lambda_2$ are the smallest and largest eigenvalues of the Hessian matrix evaluated at the target, in other words, $\nabla^2 f(p^\star,q^\star)$. Since the objective function is two-dimensional, the sum of its smallest and largest eigenvalues is equal to its trace, so
\begin{align*}
\lambda_1+\lambda_2
  = \trace(\nabla^2 f)
  = \frac{1}{n} \sum_{i=1}^n \trace(\nabla^2 f_i)
\end{align*}
where the trace of the local Hessian is
\begin{align*}
\trace(\nabla^2 f_i) = 2-\tfrac{r_i}{\sqrt{(p_i-p)^2 + (q_i-q)^2}}.
\end{align*}
The trace is equal to one at the target, so the optimal stepsize is $\alpha=2$. Since NIDS and EXTRA are unstable with this stepsize, we instead use $\alpha=1$ and $\alpha=0.5$, respectively, for these algorithms. The parameters of SVL are completely determined by $\sigma$ and $\kappa=\tfrac{1+\rho}{1-\rho}$.

We choose the contraction factor as the convergence rate of centralized gradient descent which is $\rho\approx 0.75$. Then our algorithm performs $m=6$ communication rounds per iteration. We have each agent initialize its states with its position $x_i^0=(p_i,q_i)\in\real^2$ and $y_i^0=(0,0)\in\real^2$.

In Figure~\ref{fig:target_localization_error}, we plot the error of each agent as a function of either the iteration or step. As expected from Theorem~\ref{thm:convergence_rate}, the error of our algorithm converges to zero at the same rate as centralized gradient descent~\eqref{eq:fixed_point_iteration} in terms of iterations. Our algorithm uses $m=6$ communications per iteration while NIDS, EXTRA, and SVL use only one; our algorithm is more efficient in terms of gradient evaluations, but also uses more communications than the other algorithms to obtain a solution with a given precision.

\section{Conclusion}\label{sec:conclusion}

We developed an algorithm for distributed optimization that uses the minimal amount of communication necessary such that the iterates converge to the optimizer at the same rate as centralized gradient descent in terms of the number of gradient evaluations. Furthermore, the convergence rate of our algorithm is near-optimal (in the worst-case) in terms of the number of communication rounds even though the gradient is not evaluated at each step. Such an algorithm is particularly useful when gradient evaluations are expensive relative to the cost of communication.

\bibliography{references}

\appendix

\section{Proof of Theorem~\ref{thm:convergence_rate}}\label{sec:proof_theorem1}

We now prove linear convergence of the iterates of our algorithm to the optimizer of the global objective function.

\textbf{Average and disagreement operators.} To simplify the notation, we define the \emph{average operator} $\avg : \real^{nd}\to\real^{nd}$ as
\begin{align*}
\avg(\x) \defeq (\tfrac{1}{n}\1\1^\tp\otimes I_d)\,\x
\end{align*}
along with the \emph{disagreement operator} $\dis : \real^{nd}\to\real^{nd}$ as
\begin{align*}
\dis(\x) \defeq \bigl((I_n-\tfrac{1}{n} \1\1^\tp)\otimes I_d\bigr)\,\x
\end{align*}
where $\otimes$ denotes the Kronecker product. Note that any point can be decomposed into its average and disagreement components since $\avg + \dis = I$. Also, the operators are orthogonal in that $\avg(\x)^\tp \dis(\y)=0$ for all $\x,\y\in\real^{nd}$.

\textbf{Vectorized form.} Defining the parameter $\lambda\defeq\sqrt{1-\rho^2}$, we can then write our algorithm in vectorized form as
\begin{subequations}\label{eq:alg_bf}
\begin{align}
\v^k &= \W^k(\x^k) \\
\u^k &= \v^k - \alpha\,\Df(\v^k) \\
\y^{k+1} &= \y^k + \x^k - \v^k \label{eq:alg_bf1} \\
\x^{k+1} &= \u^k - \lambda\,\y^{k+1} \label{eq:alg_bf2}
\end{align}
\end{subequations}
with $\avg(\y^0)=0$ where the concatenated vectors are
\begin{align*}
\u^k \defeq \bmat{u_1^k \\ \vdots \\ u_n^k}, \
\v^k \defeq \bmat{v_1^k \\ \vdots \\ v_n^k}, \
\x^k \defeq \bmat{x_1^k \\ \vdots \\ x_n^k}, \
\y^k \defeq \bmat{y_1^k \\ \vdots \\ y_n^k},
\end{align*}
and the $m$-step consensus operator $\W^k : \real^{nd}\to\real^{nd}$ and global gradient operator $\Df : \real^{nd}\to\real^{nd}$ are defined as\footnote{We use the over-bar in $\Df$ to distinguish it from the gradient of the global objective function $f$ in~\eqref{eq:distropt}. The operators are related by $(\tfrac{1}{n}\1^\tp\otimes I_d) \Df(\1\otimes x)=\df(x)$.}
\begin{align*}
\W^k \defeq \prod_{\ell=1}^m \bigl(W^{k,\ell}\otimes I_d\bigr)
\quad\text{and}\quad
\Df(\v) \defeq \bmat{\df_1(v_1) \\ \vdots \\ \df_n(v_n)}.
\end{align*}

\textbf{Fixed-point.} Define the points $\u^\star,\v^\star,\x^\star,\y^\star\in\real^{nd}$ as
\begin{align*}
\v^\star = \x^\star = \1\otimes x^\star, \
\u^\star = \v^\star-\alpha\Df(\v^\star), \
\y^\star = \tfrac{1}{\lambda}(\u^\star-\x^\star).
\end{align*}
Then $(\u^\star,\v^\star,\x^\star,\y^\star)$ is a fixed-point of the concatenated system~\eqref{eq:alg_bf} since the gossip matrix is doubly-stochastic at each step. Also, $\avg(\y^\star)=0$ since $x^\star$ satisfies~\eqref{eq:gradzero}.

\textbf{Error system.} To analyze the algorithm, we use a change of variables to put it in error coordinates. 
The error vectors
\begin{align*}
\bar{\u}^k &\defeq \u^k - \u^\star  &  \bar{\x}^k &\defeq \x^k - \x^\star \\
\bar{\v}^k &\defeq \v^k - \v^\star  &  \bar{\y}^k &\defeq \y^k - \y^\star
\end{align*}
satisfy the iterations
\begin{subequations}\label{eq:alg_bar}
\begin{align}
\bar{\y}^{k+1} &= \bar{\y}^k + \bar{\x}^k - \bar{\v}^k \label{eq:alg_bar1} \\
\bar{\x}^{k+1} &= \bar{\u}^k - \lambda\,\bar{\y}^{k+1} \label{eq:alg_bar2}
\end{align}
\end{subequations}
for $k\ge 0$.

\textbf{Fixed-point operator.} From Assumption~\ref{assumption:functions}, the global gradient operator $\Df$ satisfies
\begin{align*}
\avg\bigl(\Df(\x^\star)\bigr) = 0
\end{align*}
and
\begin{align}\label{eq:TT}
\|\x-\x^\star-\alpha\,(\Df(\x)-\Df(\x^\star))\| \le \rho\,\|\x-\x^\star\|
\end{align}
for all $\x\in\real^{nd}$. In other words, $I-\alpha\,\Df$ is a contraction with respect to the point $\x^\star$ with contraction factor $\rho$.

\textbf{Consensus operator.} From Assumption~\ref{assumption:W} along with the definition of $m$, the consensus operator $\W^k$ satisfies
\begin{align}\label{eq:WW}
\|\dis\bigl(\W^k(\x)\bigr)\| &\le \sigma^m\,\|\dis(\x)\| \le \sigma_0\,\|\dis(\x)\|
\end{align}
for all $\x\in\real^{nd}$ and all $k\ge 0$ where
\begin{align*}
\sigma_0\defeq\tfrac{\sqrt{1+\rho}-\sqrt{1-\rho}}{2}.
\end{align*}

\textbf{Consensus direction.} We now derive some properties of the average error vectors. Using the assumption that the gossip matrix is doubly-stochastic, we have
\begin{align}\label{eq:avg_x_v}
\avg(\bar{\x}^k) = \avg(\bar{\v}^k) \quad\forall k\ge 0.
\end{align}
The iterates are initialized such that $\avg(\bar{\y}^0)=0$ (recall that $\avg(\y^\star)=0$). Taking the average of~\eqref{eq:alg_bf1}, we have that the average is preserved. In other words, we have that $\avg(\bar{\y}^{k+1}) = \avg(\bar{\y}^k)$ for all $k\ge 0$. Then by induction,
\begin{align}\label{eq:avg_y}
\avg(\bar{\y}^k)=0 \quad\forall k\ge 0. 
\end{align}

\textbf{Lyapunov function.} To prove convergence, we will show that the function $V : \real^{nd}\times\real^{nd}\to\real$ defined by
\begin{align}\label{eq:lyap}
V(\bar{\x},\bar{\y}) \defeq \|\avg(\bar{\x})\|^2
  + \bmat{\dis(\bar{\x}) \\ \dis(\bar{\y})}^\tp\!
    \left(\bmat{1 & \lambda \\ \lambda & \lambda}\otimes I_{nd}\right)
    \bmat{\dis(\bar{\x}) \\ \dis(\bar{\y})}
\end{align}
is a Lyapunov function for the algorithm, that is, it is both positive definite and decreasing along system trajectories. Note that $\lambda\in(0,1)$ since $\rho\in(0,1)$, so the matrix in~\eqref{eq:lyap} is positive definite. Then $V$ is also positive definite, meaning that $V(\bar{\x},\bar{\y})\ge 0$ for all $\bar{\x}$ and $\bar{\y}$, and $V(\bar{\x},\bar{\y})=0$ if and only if $\bar{\x}=0$ and $\dis(\bar{\y})=0$ (recall that $\avg(\bar{\y}^k)=0$). Next, we show that the Lyapunov function decreases by a factor of at least $\rho^2$ at each iteration. Define the weighted difference in the Lyapunov function between iterations as
\begin{align*}
\Delta V^k \defeq V(\bar{\x}^{k+1},\bar{\y}^{k+1}) - \rho^2\,V(\bar{\x}^k,\bar{\y}^k).
\end{align*}
Subsituting the expressions for the iterates in~\eqref{eq:alg_bar} and using the properties of the average iterates in~\eqref{eq:avg_x_v} and~\eqref{eq:avg_y}, we have
\begin{align*}
\Delta V^k =
  &- \bigl(\rho^2\,\|\bar{\v}^k\|^2 - \|\bar{\u}^k\|^2\bigr) \nonumber \\
  &- 2\rho^2\,\bigl(\sigma_0^2\,\|\dis(\bar{\x}^k)\|^2 - \|\dis(\bar{\v}^k)\|^2\bigr) \\
  &- 2\sigma_0^2\, \bigl\|\dis\bigl(\bar{\v}^k + \lambda\,(\bar{\x}^k + \bar{\y}^k) \bigr)\bigr\|^2. \nonumber
\end{align*}
The first term is nonpositive since $\Df$ satisfies~\eqref{eq:TT}, the second since $\W^k$ satisfies~\eqref{eq:WW}, and the third since it is a squared norm. Therefore, $\Delta V^k\le 0$ for all $k\ge 0$. Applying this inequality at each iteration and summing, we obtain the bound
\begin{align*}
V(\bar{\x}^k,\bar{\y}^k) \le \rho^{2k}\,V(\bar{\x}^0,\bar{\y}^0) \quad\forall k\ge 0.
\end{align*}

\textbf{Bound.} Finally, we use the Lyapunov function to show that $\|x_i^k-x^\star\|$ converges to zero linearly with rate $\rho$ for each agent $i\in\{1,\ldots,n\}$. The norm is upper bounded by
\begin{align*}
\|x_i^k-x^\star\|^2
  &\le \cond\left(\bmat{1 & \lambda \\ \lambda & \lambda}\right)\, V(\bar{\x}^k,\bar{\y}^k)
  \le c^2\,\rho^{2k}
\end{align*}
where the nonnegative constant $c\in\real$ is defined as
\begin{align*}
c \defeq \sqrt{\cond\left(\bmat{1 & \lambda \\ \lambda & \lambda}\right)\, V(\bar{\x}^0,\bar{\y}^0)}
\end{align*}
and $\cond(\cdot)$ denotes the condition number. Taking the square root, we obtain the bound
\begin{align*}
\|x_i^k-x^\star\| \le c\,\rho^k
\end{align*}
for each agent $i\in\{1,\ldots,n\}$ and iteration $k\ge 0$. $\hfill\qed$

\end{document}